# PA is instantiationally complete, but algorithmically incomplete

## An alternative interpretation of Gödelian incompleteness under Church's Thesis that links formal logic and computability

Bhupinder Singh Anand[1]

We define instantiational and algorithmic completeness for a formal language. We show that, in the presence of Church's Thesis, an alternative interpretation of Gödelian incompleteness is that Peano Arithmetic is instantiationally complete, but algorithmically incomplete. We then postulate a Provability Thesis that links Peano Arithmetic and effective algorithmic computability, just as Church's Thesis links Recursive Arithmetic and effective instantiational computability.

## 1. Introduction: Interpreting universal quantification

An issue that, sooner or later, seems to obfuscate every philosophical discussion on the relationship between formal logic and computability, and one which seems to lie also at the root of most foundational ambiguities, is the equivocal, semantic, interpretation of universal quantification (i.e., of unqualified generalisation).

In other words, standard interpretations of classical theory do not provide an unambiguous answer to the question:

---

[1] The author is an independent scholar. E-mail: re@alixcomsi.com; anandb@vsnl.com. Postal address: 32, Agarwal House, D Road, Churchgate, Mumbai - 400 020, INDIA. Tel: +91 (22) 2281 3353. Fax: +91 (22) 2209 5091.

Can we unequivocally interpret a universally quantified formula of a language L under an interpretation M?

## 2. Differentiating instantational and algorithmic Tarskian satisfiability

We begin by noting that Tarski's definitions of the satisfiability and truth of the formulas of a formal language, say L, under an interpretation, say M, can be qualified to differentiate between at least two, essentially different, methods of assigning, for instance, the valuation TRUE to a predicate of M, by appeal to definitions such as the following[2]:

(*i*) **Effective instantiational satisfiability**: A predicate, say $R(x)$, of an interpretation M, of a language L, is effectively satisfiable instantiationally if, and only if, given any element $s$ in the domain D of the interpretation M, there is some effective[3] method, say $T(s)$, of determining that $R(s)$ is satisfied under the interpretation M;

(*ii*) **Effective algorithmic satisfiability**: A predicate, say $R(x)$, of an interpretation M, of a language L, is effectively satisfiable algorithmically if, and only if, there is some effective method, say $UT$, such that, given any element $s$ in the domain D of the interpretation M, $UT$ effectively determines that $R(s)$ is satisfied under the interpretation M.

---

[2] More generally, the definitions would be extended to cover effective instantiational, and algorithmic, methods of assigning the valuation TRUE to any predicate of M.

[3] A relation $R(x_1, x_1, \ldots, x_n)$ may be thought of as being effectively decidable if, and only if, there is a mechanical procedure - which terminates in a finite number of steps - for determining whether $R(x_1, x_1, \ldots, x_n)$ holds or not, when the arguments $x_1, x_1, \ldots, x_n$ are given.

Now, if $R(x)$ is the interpretation in M of an L-formula, $[R(x)]$[4], then (*i*) is satisfied if, for every *s* in D, $[R(s)]$ is provable[5] in L, and L is assumed to be sound[6]. Similarly, (*ii*) is satisfied if $[(Ax)R(x)]$ is provable in such an L.

Further, although every algorithmically satisfiable predicate of M is instantiationally satisfiable in M, Gödel has shown ([Go31a], p25(1)[7]) - by a constructive and intuitionistically unobjectionable meta-argument - that there are instances when $[R(x)]$ is instantiationally provable[8] if L is a consistent Peano Arithmetic, so (*i*) holds; but $[(Ax)R(x)]$ is not provable in L, so (*ii*) may not hold.

It follows that (*i*) and (*ii*) cannot be assumed equivalent without qualification.

## 3. Differentiating instantiational and algorithmic Tarskian truth

This suggests that, if we treat (*i*) as corresponding to the classical definition of Tarskian truth in M, which we may denote symbolically by $\models_M [(Ax)R(x)]$, then we can treat (*ii*) as the definition of a qualified, algorithmic, Tarskian truth in M - which we may denote symbolically by $\|=_M [(Ax)R(x)]$.

---

[4] When referring to mathematical expressions, we use square brackets to indicate that the expression within them is intended to be viewed syntactically as a formula (i.e., a string of symbols) of some formal language, say L; lack of square brackets indicates that the expression is intended to be viewed under some given interpretation, say M, of L, according to the rules governing the interpretation.

[5] We define a proof in a first-order language L as a finite sequence $[A_1]$, $[A_2]$, ... $[A_n]$ of well-formed formulas of L such that, for each *i*, either $[A_i]$ is an axiom of L, or $[A_i]$ is a direct consequence of some of the preceding well-formed formulas of the sequence by virtue of one of the rules of inference of L. We define this sequence as a proof of a well-formed formula $[R]$ of L if, and only if, $[A_n]$ is the formula $[R]$.

[6] We define a language as sound if, and only if, its provable formulas are true under every interpretation.

[7] In his Arithmetic, P, Gödel constructs a formula, say $[(Ax)R(x, p)]$, with Gödel-number 17Gen*r*, that is not provable in a consistent P, although $[R(x, p)]$, whose Gödel-number is *r*, is P-provable for every substitution of the variable $[x]$ by any numeral, $[n]$, of P.

[8] By '$[R(x)]$ is instantiationally provable' we mean that the L-formula, $[R(n)]$, obtained by substituting the numeral $[n]$ for the variable $[x]$, is L-provable for every numeral in L.

So we can, without loss of generality, symbolically express the distinction between instantiational Tarskian truth, and algorithmic Tarskian truth, by defining:

(*1*)    $\models_M [(Ax)R(x)] \Longleftrightarrow$ (*i*);

and:

(*2*)    $\Vert =_M [(Ax)R(x)] \Longleftrightarrow$ (*ii*);

where:

(*3*)    $\Vert =_M [(Ax)R(x)] \Longrightarrow \models_M [(Ax)R(x)]$;

but the following does not always hold:

(*4*)    $\models_M [(Ax)R(x)] \Longrightarrow \Vert =_M [(Ax)R(x)]$.

## 4. Defining instantiational and algorithmic provability

We can, then, define:

(*iii*) **Instantiational provability**: An instantiationally Tarskian-true predicate, $R(\ldots)$[9], of an interpretation, say M, of a formal language, say L, is instantiationally provable in L if, and only if, $R(\ldots)$ is the interpretation of an L-formula, $[R(\ldots)]$, and every instantiation of $R(\ldots)$ in M is the interpretation of an L-provable formula.

(*iv*) **Algorithmic provability**: An algorithmically Tarskian-true predicate, $R(\ldots)$, of an interpretation, say M, of a formal language, say L, is algorithmically provable in L if, and only if, $R(\ldots)$ is the interpretation of an L-formula, $[R(\ldots)]$, and $[(A\ldots)R(\ldots)]$ is provable in L.

---

[9] In the following, we denote an expression of $n$-arguments, say $R(x_1, x_2, \ldots, x_n)$, by $R(\ldots)$, and its closure under universal quantification, $(A x_1)(A x_2)\ldots(A x_n)R(x_1, x_1, \ldots, x_1)$, by $(A\ldots)R(\ldots)$.

## 5. Defining instantiational and algorithmic completeness

We can also define, further:

**Instantiational completeness**: A language L is instantiationally complete with respect to an interpretation M if, and only if, every instantiationally Tarskian-true predicate, $R(\ldots)$, of M is instantiationally equivalent to a predicate, $R'(\ldots)$, of M, every instantiation of which is the interpretation in M of an L-provable formula.

**Algorithmic completeness**: A language L is algorithmically complete with respect to an interpretation M if, and only if, every algorithmically Tarskian-true predicate, $R(\ldots)$, of M is provable in L.

We now show the significance of the above definitions, by considering the particular case where L is a consistent first-order Peano Arithmetic, and M its standard interpretation.

## 6. Differentiating between instantiational and algorithmic effective computability

Now, a standard expression of Church's Thesis is:

**Church's Thesis**: A (partial) number-theoretic function is effectively computable if, and only if, it is (partial) recursive.

If we, again, differentiate between instantiational and algorithmic computability, we can interpret the above as asserting that:

**Instantiational Church's Thesis**: A (partial) number-theoretic function (or relation, treated as a Boolean function) is effectively computable instantiationally if, and only if, it is instantiationally equivalent to a (partial) recursive function (or relation, treated as a Boolean function).

and:

> **Algorithmic Church's Thesis**: A (partial) number-theoretic function (or relation, treated as a Boolean function) is effectively computable algorithmically if, and only if, it is (partial) recursive.

## 7. Is first-order Peano Arithmetic instantiationally complete with respect to its standard interpretation?

Now, by the Instantiational Church's Thesis, every instantiationally Tarskian-true predicate, $R(\ldots)$, of M is instantiationally equivalent to a recursive relation, $R'(\ldots)$, of M that, treated as a Boolean function, always evaluates as TRUE.

By Gödel's Theorems V and VII [Go31a][10], it follows that there is, thus, a PA-formula, $[R''(\ldots)]$, such that its standard interpretation in M, the arithmetical relation $R''(\ldots)$, is instantiationally equivalent, in M, to $R'(\ldots)$, and, hence, to $R(\ldots)$; and every instantiation of $R''(\ldots)$ in M is the standard interpretation of a PA-provable formula.

It follows that:

> **Meta-theorem 1**: Church's Thesis implies that first-order Peano Arithmetic is instantiationally complete with respect to its standard interpretation.

Thus, Church's Thesis is the meta-postulate that provides the arithmetic consequences, in first-order Peano Arithmetic, which are intended by the second-order Induction Axiom of Peano's Postulates ([Me64], p102).

---

[10] We also refer to the combination of these two theorems as Hilbert and Bernays Representation Theorem.

The Thesis, essentially, postulates that any arbitrary property that holds instantially under the standard interpretation of Peano's second-order Induction Axiom is instantiationally equivalent to a recursive relation in the interpretation.

It follows that:

> **Corollary 1.1**: Church's Thesis implies that first-order Peano Arithmetic formalizes Dedekind's Peano Postulates completely with respect to Tarskian truth.

The question arises: How do we interpret Gödelian incompleteness in the presence of Church's Thesis?

## 8. Is Peano Arithmetic algorithmically complete with respect to its standard interpretation?

Now, Peano Arithmetic would, further, be algorithmically complete with respect to its standard interpretation, M, if, and only if, every algorithmically Tarskian-true predicate, $R(\ldots)$, of M were provable in PA.

In the Appendix we show, however, that this is not the case, since Gödel's Tarskian-true primitive recursive predicate, $\sim xB(Sb(p\ 19|Z(p)))$[11] - which is algorithmic by definition - is not the standard interpretation of any of its representations in Gödel's formal system P (and the argument remains valid for any first-order Peano Arithmetic).

We thus have:

> **Meta-theorem 2**: First-order Peano Arithmetic is not algorithmically complete with respect to its standard interpretation, M.

---

[11] This is the primitive recursive relation whose formal representation in Gödel's Arithmetic, P, we denote by $[R(x, p)]$. We define these expressions fully in the Appendix.

So, in the presence of Church's Thesis, the Gödelian incompleteness of a Peano Arithmetic is, essentially, the assertion that the Arithmetic is not algorithmically complete with respect to its standard interpretation, M.

## 9. Is Peano Arithmetic algorithmically complete with respect to the Arithmetical part of M?

We define:

> **Arithmetic completeness**: An Arithmetic, A, is arithmetically complete with respect to its standard interpretation, M, if, and only if, every algorithmically Tarskian-true arithmetical predicate, $R(\ldots)$, of M is provable in A.

The question arises: Is PA arithmetically complete?

Now, it seems reasonable, and intuitively unobjectionable, to assume that, if there is an algorithm such that every instantiation of an arithmetical predicate, $R(\ldots)$, is effectively decidable as Tarskian-true in M, then there must be a corresponding PA-proof sequence for the formula, $[R(\ldots)]$.

We may formally define this assertion as a:

> **Provability Thesis**: Every algorithmically Tarskian-true arithmetical predicate, $R(\ldots)$, is provable in PA.

It follows that:

> **Meta-theorem 3**: The Provability Thesis implies that first-order Peano Arithmetic is arithmetically complete with respect to its standard interpretation, M.

## 10. Is the Provability Thesis effectively verifiable?

The question arises: Is the Provability Thesis effectively verifiable?

Now, if it were always effectively decidable in M whether or not an arbitrary number-theoretic relation is algorithmically Tarskian-true - i.e., it is effectively satisfied algorithmically in the sense of (*ii*) - then, assuming Church's Thesis, it would follow that Turing's Halting Problem would be decidable by a Turing machine.

Since this is not the case, it follows that the Provability Thesis, like Church's Thesis, may not be effectively verifiable.

## 11. Soundness

Now, standard interpretations of classical logic define the following semantic meta-implication as an instance of the property of the 'soundness' of PA:

(*5*)        $\vdash_{PA} [(Ax)R(x)] \Rightarrow \models_{M} [(Ax)R(x)]$.

Further, standard interpretations of classical theory show that the converse meta-implication, $\models_{M} [(Ax)R(x)] \Rightarrow \vdash_{PA} [(Ax)R(x)]$, does not necessarily hold - a consequence of Tarski's Theorem that arithmetical truth cannot be defined arithmetically.

The Provability Thesis is, thus, an attempt to offer a possible answer to the question: Can, and if so, when, does the converse hold?

It is the postulation that, if $R(x)$ is an arithmetical predicate, then, in the sense of (*2*):

(*6*)        $\models_{M} [(Ax)R(x)] \Rightarrow \vdash_{PA} [(Ax)R(x)]$

 The detailed consequences of such a postulation are beyond the immediate intent, and scope, of the subject of this essay.

## 12. Some consequences of postulating a Provability Thesis for PA

However, we note, briefly, that, if we postulate a Provability Thesis for Peano Arithmetic, then:

(*a*) we would interpret Gödel's reasoning as showing that, whereas (*ii*) implies (*i*), the converse does not always hold;

(*b*) Peano Arithmetic would be ω-inconsistent[12] - and, so, Gödel's Theorem VI [Go31a], and its consequences, would hold vacuously[13];

> According to the Provability Thesis, there is no algorithm for determining that, given any natural number $n$ in the standard interpretation M of PA, $R(n)$ holds in M, although $R(n)$ is true in M for every natural number $n$, and $[R(n)]$ is provable in PA for every numeral $[n]$.

> Since mere addition of $[(Ax)R(x)]$ as an axiom to PA would not entail the introduction of an algorithm in M for determining that, given any natural number $n$, $R(n)$ holds in M, it follows that PA+$[(Ax)R(x)]$ would be inconsistent under the standard interpretation.

> It follows that, if PA is consistent, then $[(Ax)R(x)]$ would not be independent of the axioms of PA. Hence, $[\sim(Ax)R(x)]$ would be PA-provable[14].

---

[12] A language L is said to be ω-inconsistent if, and only if, for every formula $[R(x)]$ of L, if $[R(n)]$ is provable for every numeral $[n]$, then it is not the case that $[(Ex)\sim R(x)]$ is also L-provable.

[13] In Theorem VI of his seminal 1931 paper, Gödel argues that, if his Arithmetic P is assumed consistent, then $[(Ax)R(x, p)]$ is not P-provable; whilst, if P is, further, assumed also ω-inconsistent, then $[\sim(Ax)R(x, p)]$, too, is not P-provable - where these expressions are as defined in the Appendix.

[14] By definition, $[\sim(Ax)R(x)]$ interprets in M as the negation of the interpretation of $[(Ax)R(x)]$ in M. Hence the PA-provability of $[\sim(Ax)R(x)]$ is the meta-assertion that there is no algorithm in M for determining that, given any natural number $n$, $R(n)$ holds in M.

Since [$R(x)$] is provable in PA whenever we substitute any numeral [$n$] for the variable [$x$], it would then follow that PA would be consistent, but not ω-consistent.

(*c*) the intuitionist objection, to concluding that the provability of [(E$x$)$R'(x)$] in PA always implies the existence of some $s$, in M, such that $R'(s)$ holds in M, would be vindicated if PA is ω-consistent; and Rosser's proof of undecidability would fail.

Since [(E$x$)$R'(x)$] is simply an abbreviation of [~(A$x$)~$R'(x)$], under the Provability Thesis, the PA provability of [~(A$x$)~$R'(x)$] is the meta-assertion that there is no algorithm in M for determining that, given any natural number $n$, ~$R'(n)$ holds in M.

However, if PA is ω-consistent, we may have that [~(A$x$)~$R'(x)$] is PA-provable, and also that [~$R(x)$] is provable in PA whenever we substitute any numeral [$n$] for the variable [$x$].

Hence we may not conclude from the PA provability of [~(A$x$)~$R'(x)$] that there is always some $s$, in M, such that $R'(s)$ holds in M.

Now Rosser's proof - that that his proposition, say [(A$x$)$R*(x)$], can be shown to be undecidable in PA without assuming ω-consistency - assumes that, from the premise [~(A$x$)$R*(x)$] in PA, we may conclude the existence of some natural number $s$ in M such that ~$R*(s)$ holds in M. However, under the Provability Thesis, such an assumption may be invalid.

(*d*) the Turing Thesis would no longer hold, since, if [(A$x$)$G(x)$] represents Gödel's undecidable proposition, then the arithmetical predicate $G(n)$, treated as a Boolean function, would be an instantiationally (assumed effectively) computable function that would not, however, be Turing-computable;

(*e*) consequently, $G(n)$ would be in the complexity class NP, but not in P;

(*f*) Hilbert's proposed ω-rule could, then, be viewed as also an attempt at a complete, semantic, definition of provability for Peano Arithmetic.

## 13. Significance of a Provability Thesis

Now, the Provability Thesis, is, essentially, the assertion that a PA formula is provable if, and only if, any one of its interpretations is Tarskian-decidable algorithmically, in the sense of (*ii*) above.

The significance of such a link between formal provability and computability emerges when it is seen in a broader perspective, where we note that:

(*a*) Church's Thesis, in its instantiational, weakened, form, (i.e., as an equivalence which is implied by the standard form of Church's Thesis, although the converse does not hold) builds an iff bridge between instantiational (assumed effective) computability and Recursive Arithmetic;

eg., A number-theoretic predicate $R(x)$, is instantiationally (assumed effectively) decidable for every natural number value of $x$ if, and only if, $R(x)$ is instantiationally (assumed effectively) equivalent to a recursive predicate $Q(x)$.

(*b*) Hilbert & Bernays Representation Theorem, and its converse, together, build an iff bridge between Recursive Arithmetic and Peano Arithmetic;

eg., A number-theoretic predicate $Q(x)$ is recursive if, and only if, $Q(x)$ is representable in PA when treated as a Boolean function.

(*c*) The Church-Turing Theorem builds an iff bridge between Recursive Arithmetic and Turing computability;

eg., A number-theoretic predicate $Q(x)$ is recursive if, and only if, $Q(x)$, treated as a Boolean function, is total and Turing-computable.

(*d*) Markov algorithms build an iff bridge between Turing computability and algorithmic computability.

eg., A number-theoretic predicate $Q(x)$, treated as a Boolean function, is Turing-computable if, and only if, Q(x), treated as a Boolean function, is partially Markov-computable.

However:

(*e*) Soundness merely builds an only-if bridge between Peano Arithmetic and Turing (algorithmic) computability;

eg., $[R(x)]$ is a PA-provable formula only if, under the standard interpretation of PA, the interpreted arithmetical predicate, $R(n)$, treated as a Boolean function, is Turing-computable as TRUE for any given natural number $n$.

whilst:

(*f*) The Provability Thesis complements, and completes, the concept of Soundness in arithmetic by building an iff bridge between Peano Arithmetic and Turing (algorithmic) computability.

eg., $[R(x)]$ is a PA-provable formula if, and only if, under the standard interpretation of PA, the interpreted arithmetical predicate, $R(n)$, treated as a Boolean function, is Turing-computable as TRUE for any given natural number $n$.

## 14. Conclusion

We have shown that Church's Thesis implies that a first-order Peano Arithmetic is instantiationally complete, in the sense that it completely formalizes Dedekind's Peano Postulates with respect to Tarskian-truth.

We have shown, further, that, in the presence of Church's Thesis, the Gödelian incompleteness of a Peano Arithmetic is, essentially, the assertion that the Arithmetic is not algorithmically complete with respect to its standard interpretation, M.

We have, then, defined a Provability Thesis that complements, and completes, the concept of Soundness in arithmetic by building an iff bridge between Peano Arithmetic and Turing (algorithmic) computability.

We have also briefly indicated how the consequences of such a Thesis - whose detailed consideration lies beyond the intent and scope of the present investigation - may be significant for standard interpretations of classical theory.

## Appendix

### Notation

Unless specified otherwise, we generally follow the notation introduced by Mendelson in his English translation of Gödel's 1931 paper [Go31a]; however, for convenience of exposition, we refer to it as Gödel's notation. Two notable exceptions: we use the notation "$(Ax)$", whose classical, standard, interpretation is "for all $x$", to denote Gödel's special symbol for Generalisation; the successor symbol is denoted by "$S$", instead of by "$f$".

Following Gödel (cf. [Go31a], footnote 13), we use square brackets to indicate that the expression $[(Ax)]$, including square brackets, only denotes the uninterpreted string[15] named[16] within the brackets. Thus, $[(Ax)]$ is not part of the formal system P, and would be replaced by Gödel's special symbolism for Generalisation wherever it occurs.

Following Gödel's definitions of well-formed formulas[17], we note that juxtaposing the string $[(Ax)]$ and the formula[18] $[F(x)]$ is the formula $[(Ax)F(x)]$, juxtaposing the symbol $[\sim]$ and the formula $[F]$ is the formula $[\sim F]$, and juxtaposing the symbol $[v]$ between the formulas $[F]$ and $[G]$ is the formula $[FvG]$.

The number-theoretic functions and relations in the following are defined explicitly by Gödel [Go31a]. The formulas are defined implicitly by his reasoning.

**Preliminary Definitions**

We take P to be Gödel's formal system[19], and define ([Go31a], Theorem VI, p24-25):

(*i*)    $Q(x, y)$ as Gödel's recursive relation $\sim xB(Sb(y\ 19|Z(y)))$[20].

---

[15] We define a "string" as any concatenation of a finite set of the primitive symbols of the formal system under consideration.

[16] We note that the "name" inside the square brackets only serves as an abbreviation for some string in P.

[17] We note that all well-formed formulas of P are strings of P, but all strings of P are not well-formed formulas of P.

[18] By "formula", we mean a "well-formed formula" as defined by Gödel ([Go31a], p11).

[19] Gödel ([Go31a], p9-13).

[20] We follow Gödel's definition of recursive number-theoretic functions and relations ([Go31a], p14-17). We note, in particular, that Gödel's recursive number-theoretic function $Sb(x\ 19|Z(y))$ is defined as the Gödel-number of the P-formula that is obtained from the P-formula whose Gödel-number is $x$ by substituting the numeral $[y]$, whose Gödel-number is $Z(y)$, for the variable whose Gödel-number is 19 wherever the latter occurs free in the P-formula whose Gödel-number is $x$ ([Go31a], p20, def.31). We also note that Gödel's recursive number-theoretic relation $xBy$ holds if, and only if, $x$ is the Gödel-number of a proof sequence for the P-formula whose Gödel-number is $y$ ([Go31a], p22, def. 45).

(*ii*)　[*R*(*x*, *y*)] as a formula that represents *Q*(*x*, *y*) in the formal system P.

(*iii*)　*q* as the Gödel-number[21] of the formula [*R*(*x*, *y*)] of P.

(*iv*)　*p* as the Gödel-number of the formula [(A*x*)*R*(*x*, *y*)][22] of P.

(*v*)　[*p*] as the numeral that represents the natural number *p* in P.

(*vi*)　*r* as the Gödel-number of the formula [*R*(*x*, *p*)] of P.

(*vii*)　17*Genr* as the Gödel-number of the formula [(A*x*)*R*(*x*, *p*)] of P.

(*viii*)　*Neg*(17*Genr*)[23] as the Gödel-number of the formula [~(A*x*)*R*(*x*, *p*)] of P.

(*ix*)　*R*(*x*, *y*) as the standard interpretation of the formula [*R*(*x*, *y*)] of P.

**Meta-theorem 2**[24]: There is a recursive relation that is not the standard representation of any of its representations in P.

*Proof*: We consider Gödel's primitive recursive relation ~*xB*(*Sb*(*y* 19|*Z*(*y*))).

(*a*)　We assume that every recursive number-theoretic function or relation is the standard interpretation of at least one of its formal representations in P[25].

---

[21] By the "Gödel-number" of a formula of P, we mean the natural number corresponding to the formula in the 1-1 correspondence defined by Gödel ([Go31a], p13).

[22] We note that "[(A*x*)][*R*(*x*, *y*)]" and "[(A*x*)*R*(*x*, *y*)]" denote the same formula of P.

[23] We note that Gödel's recursive number-theoretic function *Neg*(*x*) is the Gödel-number of the P-formula that is the negation of the P-formula whose Gödel-number is *x* ([Go31a], p18, def. 13).

[24] Meta-theorem 2 is, essentially, the meta-thesis that, if we assume the primitive recursive relation ~*xB*(*Sb*(*y* 19|*Z*(*y*))) to be an abbreviation of some formula of P, then we arrive at an inconsistency in P.

[25] In other words, we assume that, if we use Gödel's recursive definitions ([Go31a], p17-20), and follow the reasoning he outlines in Theorem V ([Go31a], p23), we can transform the relation ~*xB*(*Sb*(*y* 19|*Z*(*y*))) into a form such that all the symbols that occur in it are standard interpretations of primitive symbols of P. See also Gödel's remarks ([Go31a], p11, footnote 22) in this context.

(*b*)  There is, thus, some P-formula [~*xB*(*Sb*(*y* 19|*Z*(*y*)))] whose standard interpretation is the primitive recursive relation ~*xB*(*Sb*(*y* 19|*Z*(*y*)))[26].

(*c*)  Now, in every model M[27] (cf. [Me64], p192-3) of P, we can also interpret[28]:

   (*i*)    the integer 0 as the interpretation of the symbol "0";

   (*ii*)   the successor operation as the interpretation of the successor function "'";

   (*iii*)  addition and multiplication as the interpretations of "+" and ".";

   (*iv*)  the interpretation of the predicate letter = as the identity relation.

(*d*)  Since the numerals of P interpret as a sub-domain of every model M of P, the natural numbers are, then, a sub-domain of every M.

(*e*)  Further, by the hypothesis (*a*), all of Gödel's 45 primitive recursive functions and relations ([Go31a], p17-22) are also, then, mirrored in every model M of P, and the P-formula, [~*xB*(*Sb*(*p* 19|*Z*(*p*)))], always interprets as the M-relation ~*xB*(*Sb*(*p* 19|*Z*(*p*))), where [*p*] is the numeral that represents the natural number *p* in P, and *p* is the Gödel-number of [(A*x*)~*xB*(*Sb*(*y* 19|*Z*(*y*)))].

(*f*)  Hence, in every model M of P, the relation ~*xB*(*Sb*(*p* 19|*Z*(*p*))) holds in M if, and only if, *x* is a M-number that is not the Gödel-number of a proof of [(A*x*)~*xB*(*Sb*(*p* 19|*Z*(*p*)))] in P.

---

[26] Since, by the hypothesis (*a*), every recursive number-theoretic function and relation is "mirrored" in P.

[27] We follow Mendelson's definition of a model ([Me64], p51): An interpretation is said to be a model for a set *T* of well-formed formulas of P if, and only if, every well-formed formula in *T* is true for the interpretation.

[28] Cf. Mendelson ([Me64], p107).

(*g*) Further, since the Gödel-number of a proof of [(A*x*)~*xB*(*Sb*(*p* 19|*Z*(*p*)))] in P is necessarily a natural number, ~*xB*(*Sb*(*p* 19|*Z*(*p*))) holds in every model M of P if *x* is an M-number that is not a natural number.

(*h*) Now, Gödel has shown (cf. [Go31a], Theorem VI) that ~*xB*(*Sb*(*p* 19|*Z*(*p*))) holds over the domain of the natural numbers.

(*i*) It follows that ~*xB*(*Sb*(*p* 19|*Z*(*p*))) is satisfied by all *x* in every model M of P.

(*j*) Hence, the P-formula [(A*x*)~*xB*(*Sb*(*p* 19|*Z*(*p*)))] is true (cf. [Me64], p51) in every model M of P, and, by Gödel's Completeness Theorem ([Me64], Corollary 2.14, p68), [(A*x*)~*xB*(*Sb*(*p* 19|*Z*(*p*)))] is P-provable.

However, since Gödel has shown that [(A*x*)~*xB*(*Sb*(*p* 19|*Z*(*p*)))] is not P-provable (cf. [Go31a], Theorem VI), we conclude that assumption (*a*) does not hold. This proves the theorem.

(*Created: Saturday*, 2*nd July 2005*, 8:16:55 *PM IST. Updated Monday*, 4*th July 2005*, 12:44:06 *AM IST, by re@alixcomsi.com.*)